\documentclass[12pt,thmsa]{article}
\usepackage{amsmath}
\usepackage{amssymb}
\usepackage{amsfonts}
\usepackage{graphicx}
\usepackage[T1]{fontenc}

\input{tcilatex}

\begin{document}

\title{Derivative Operators in Metric and Geometric Structures}
\author{{\footnotesize V. V. Fern\'{a}ndez}$^{1}${\footnotesize , A. M. Moya}%
$^{2}${\footnotesize , W. A. Rodrigues Jr}.$^{1}${\footnotesize \ } \\
$^{1}\hspace{-0.1cm}${\footnotesize Institute of Mathematics, Statistics and
Scientific Computation}\\
{\footnotesize \ IMECC-UNICAMP CP 6065}\\
{\footnotesize \ 13083-970 Campinas, SP, Brazil }\\
{\footnotesize e-mail:walrod@ime.unicamp.br ; virginvelfe@accessplus.com.ar}%
\\
{\footnotesize \ }$^{2}${\footnotesize Department of Mathematics, University
of Antofagasta, Antofagasta, Chile} \\
{\footnotesize e-mail: mmoya@uantof.cl} }
\maketitle

\begin{abstract}
This paper (the seventh paper in a series of eight) continues the
development of our theory of multivector and extensor calculus on smooth
manifolds. Here we deal first with the concepts of ordinary Hodge
coderivatives, duality identities, and Hodge coderivative identities. Then,
we recall the concept of a Levi-Civita geometric structure and the concepts
of Levi-Civita and gauge derivatives. New formulas that are important in the
Lagrangian theory of multivector and extensor fields are obtained. We
introduce also he concept of \emph{covariant} Hodge coderivative. We detail
how all these concepts are related.
\end{abstract}

\tableofcontents

\section{Introduction}

This is the seventh paper in a series of eight. Here, we continue the
development of our theory of multivector and extensor calculus on smooth
manifolds. We introduce in Section 2 the ordinary Hodge coderivatives,
duality identities, and Hodge coderivative identities. In Section 3 we
recall the concept of a Levi-Civita geometric structure and the concepts of
Levi-Civita and \textit{gauge }derivatives. Several important formulas that
appear in the Lagrangian formulation of the theory of multivector and
extensor fields on smooth manifolds are obtained.\footnote{%
In \cite{1} we gave a preliminary presentation of the Lagrangian theory of
multivector and extensor fields on Minkowski spacetime. A more general
theory based on the developments of the present series of papers is in
preparation.} In Section 4 we introduce the new concept of \emph{covariant}
Hodge coderivative. We study in details how all these important concepts are
related. Finally in Section 5 we present our conclusions.

\section{Ordinary Hodge Coderivatives}

Let $U$ be an open subset of $U_{o},$ and let $(U,g)$ be a \emph{metric
structure} on $U.$ Let us take any pair of \emph{reciprocal frame fields} on 
$U,$ say $(\{e_{\mu}\},\{e^{\mu}\}),$ i.e., $e_{\mu}\cdot
e^{\nu}=\delta_{\mu }^{\nu}.$ In particular, the \emph{fiducial frame field}%
, namely $\{b_{\mu }\},$ due to its \emph{orthonormality}, i.e., $%
b_{\mu}\cdot b_{\nu}=\delta_{\mu\nu},$ is a \emph{self-reciprocal frame field%
}, i.e., $b^{\mu }=b_{\mu}.$ We note also that $\{b_{\mu}\}$ is an \emph{%
ordinarily constant frame field} on $U,$ i.e., 
\begin{equation}
a\cdot\partial_{o}b_{\mu}=0,\text{ for each }\mu=1,\ldots,n.  \label{OHD.1}
\end{equation}

Associated to $(\{e_{\mu}\},\{e^{\mu}\}),$ the smooth pseudoscalar field on $%
U,$ namely $\tau,$ defined by 
\begin{equation}
\tau=\sqrt{e_{\wedge}\cdot e_{\wedge}}e^{\wedge},  \label{OHD.2}
\end{equation}
where $e_{\wedge}=e_{1}\wedge\ldots\wedge e_{n}\in\mathcal{M}^{n}(U)$ and $%
e^{\wedge}=e^{1}\wedge\ldots\wedge e^{n}\in\mathcal{M}^{n}(U),$ is said to
be the \emph{standard volume pseudoscalar field }for the \emph{local
coordinate system} $(U_{o},\phi_{o}).$

Such $\tau\in\mathcal{M}^{n}(U)$ has the fundamental property 
\begin{equation}
\tau\cdot\tau=\tau\lrcorner\widetilde{\tau}=\tau\widetilde{\tau}=1.
\label{OHD.2a}
\end{equation}
It follows from the obvious result $e_{\wedge}\cdot e^{\wedge}=1.$

>From Eq.(\ref{OHD.2a}) we can get an expansion formula for smooth
pseudoscalar fields on $U,$ i.e., 
\begin{equation}
I=(I\cdot\tau)\tau.  \label{OHD.2b}
\end{equation}

In particular, because of the obvious properties $b_{\wedge}\cdot b_{\wedge
}=1$ and $b^{\wedge}=b_{\wedge},$ the standard volume pseudoscalar field
associated to $\{b_{\mu}\}$ is just $b_{\wedge}.$ It will be called the 
\emph{canonical volume pseudoscalar field} for $(U_{o},\phi_{o}).$

We emphasize that $b_{\wedge}$ is an \emph{ordinarily constant }smooth
pseudoscalar field on $U,$ i.e., 
\begin{equation}
a\cdot\partial_{o}b_{\wedge}=0.  \label{OHD.3}
\end{equation}
Eq.(\ref{OHD.3}) can be proved by using Eq.(\ref{OHD.1}) and the Leibnitz
rule for the exterior product of smooth multivector fields.

By using Eq.(\ref{OHD.3}) and the general Leibnitz rule for $a\cdot
\partial_{o}$ we can deduce the remarkable property 
\begin{equation}
a\cdot\partial_{o}(b_{\wedge}\ast X)=b_{\wedge}\ast(a\cdot\partial_{o}X),
\label{OHD.4}
\end{equation}
where $\ast$ means either exterior product or any canonical product of
smooth multivector fields.

On the other hand we have that all $(\{e_{\mu}\},\{e^{\mu}\})$ must be
necessarily an \emph{extensor-deformation} of $\{b_{\mu}\}$. This statement
means that there exists a non-singular smooth $(1,1)$-extensor field on $U,$
say $\varepsilon,$ such that 
\begin{align}
e_{\mu} & =\varepsilon(b_{\mu}),  \label{OHD.5a} \\
e^{\mu} & =\varepsilon^{*}(b_{\mu}),\text{ for each }\mu=1,\ldots,n.
\label{OHD.5b}
\end{align}

Then, by putting Eq.(\ref{OHD.5a}) and Eq.(\ref{OHD.5b}) into Eq.(\ref{OHD.2}%
), we have that 
\begin{equation*}
\tau=(\underline{\varepsilon}(b_{\wedge})\cdot\underline{\varepsilon }%
(b_{\wedge}))^{1/2}\underline{\varepsilon}^{\ast}(b_{\wedge})=(\left.
\det\right. ^{2}[\varepsilon]b_{\wedge}\cdot b_{\wedge})^{1/2}\left.
\det\right. ^{-1}[\varepsilon]b_{\wedge}=sgn(\det[\varepsilon])b_{\wedge},
\end{equation*}
i.e., 
\begin{equation}
\tau=\pm b_{\wedge}.  \label{OHD.6}
\end{equation}

>From Eq.(\ref{OHD.3}) and Eq.(\ref{OHD.4}), by taking into account Eq.(\ref%
{OHD.6}), we have two remarkable properties 
\begin{align}
a\cdot\partial_{o}\tau & =0,  \label{OHD.6a} \\
a\cdot\partial_{o}(\tau\ast X) & =\tau\ast(a\cdot\partial_{o}X),\text{ for
all }X\in\mathcal{M}(U).  \label{OHD.6b}
\end{align}

Associated to $(\{e_{\mu}\},\{e^{\mu}\}),$ the smooth pseudoscalar field on $%
U,$ namely $\underset{g}{\tau},$ defined by 
\begin{equation}
\underset{g}{\tau}=\sqrt{\left| e_{\wedge}\underset{g}{\cdot}e_{\wedge
}\right| }e^{\wedge}=\sqrt{\left| \det[g]\right| }\tau,  \label{OHD.7}
\end{equation}
will be said to be a \emph{metric volume pseudoscalar field} for $%
(U_{o},\phi_{o}).$

Such $\underset{g}{\tau}$ $\in\mathcal{M}^{n}(U)$ satisfies the basic
property 
\begin{equation}
\underset{g}{\tau}\underset{g^{-1}}{\cdot}\underset{g}{\tau}=\tau \underset{%
g^{-1}}{\lrcorner}\underset{g}{\widetilde{\tau}}=\underset{g}{\tau }\underset%
{g^{-1}}{}\underset{g}{\widetilde{\tau}}=(-1)^{q}.  \label{OHD.7a}
\end{equation}
In order to prove it we should recall that $sgn(\det[g])=(-1)^{q},$ where $q$
is the \emph{number of negative eigenvalues} of $g.$

An expansion formula for smooth pseudoscalar fields on $U$ can be also
obtained from Eq.(\ref{OHD.7a}), i.e., 
\begin{equation}
I=(-1)^{q}(I\underset{g^{-1}}{\cdot}\underset{g}{\tau})\underset{g}{\tau}.
\label{OHD.7b}
\end{equation}

Associated to $\tau,$ the smooth extensor field on $U,$ namely $\star,$
defined by 
\begin{equation}
\star X=\widetilde{X}\lrcorner\tau=\widetilde{X}\tau,  \label{OHD.8}
\end{equation}
will be called the \emph{standard Hodge extensor field }on $U.$

Such $\star$ is \emph{non-singular} and its inverse $\star^{-1}$ is given by 
\begin{equation}
\star^{-1}X=\tau\llcorner\widetilde{X}=\tau\widetilde{X}.  \label{OHD.8a}
\end{equation}

By using a property analogous to that result given by Eq.(\ref{OHD.6b}) we
can easily prove that the standard Hodge extensor field is \emph{ordinarily
constant}, i.e., 
\begin{equation}
a\cdot\partial_{o}\star=0.  \label{OHD.8b}
\end{equation}
Also, 
\begin{equation}
a\cdot\partial_{o}\star^{-1}=0.  \label{OHD.8c}
\end{equation}

Associated to $\underset{g}{\tau},$ the smooth extensor field on $U,$ namely 
$\underset{g}{\star},$ defined by 
\begin{equation}
\underset{g}{\star}X=\widetilde{X}\underset{g^{-1}}{\lrcorner}\underset{g}{%
\tau}=\widetilde{X}\underset{g^{-1}}{}\underset{g}{\tau}=\sqrt{\left| \det[g]%
\right| }\underline{g}^{-1}(\widetilde{X})\lrcorner\tau=\sqrt{\left| \det[g]%
\right| }\underline{g}^{-1}(\widetilde{X})\tau,  \label{OHD.9}
\end{equation}
will be called the \emph{metric Hodge extensor field }on $U.$ Of course, it
is associated to the metric structure $(U,g).$

Such $\underset{g}{\star}$ is also non-singular and its inverse, namely $%
\underset{g}{\star}^{-1},$ is given by 
\begin{align}
\underset{g}{\star}^{-1}X & =(-1)^{q}\underset{g}{\tau}\underset{g^{-1}}{%
\llcorner}\widetilde{X}=(-1)^{q}\underset{g}{\tau}\underset{g^{-1}}{}%
\widetilde{X},  \notag \\
& =(-1)^{q}\sqrt{\left| \det[g]\right| }\tau\llcorner\underline{g}^{-1}(%
\widetilde{X})=(-1)^{q}\sqrt{\left| \det[g]\right| }\tau\underline {g}^{-1}(%
\widetilde{X}).  \label{OHD.9a}
\end{align}

\subsection{Duality Identities}

We present in this subsection two interesting and useful formulas which
relate the \emph{ordinary curl} $\partial_{o}\wedge$ to the \emph{ordinary
contracted divergence} $\partial_{o}\lrcorner.\vspace{0.1in}$

\textbf{i.} For all $X\in\mathcal{M}(U)$ it holds 
\begin{equation}
\tau(\partial_{o}\wedge X)=(-1)^{n+1}\partial_{o}\lrcorner(\tau X).
\label{DI.1}
\end{equation}

\textbf{Proof}

We will use the so-called duality identity $I(a\wedge
Y)=(-1)^{n+1}a\lrcorner(IY),$ where $a\in\mathcal{U}_{o},$ $I\in\bigwedge^{n}%
\mathcal{U}_{o}$ and $Y\in\bigwedge\mathcal{U}_{o}.$ By recalling the known
identities $\partial_{a}\wedge(a\cdot\partial_{o}X)=\partial_{o}\wedge X$
and $\partial_{a}\lrcorner(a\cdot\partial_{o}X)=\partial_{o}\lrcorner X,$
and using Eq.(\ref{OHD.6b}), we get 
\begin{align*}
\tau(\partial_{o}\wedge X) & =\tau(\partial_{a}\wedge(a\cdot\partial
_{o}X))=(-1)^{n+1}\partial_{a}\lrcorner(\tau(a\cdot\partial_{o}X)), \\
& =(-1)^{n+1}\partial_{a}\lrcorner(a\cdot\partial_{o}(\tau
X))=(-1)^{n+1}\partial_{o}\lrcorner(\tau X).\blacksquare
\end{align*}

\textbf{ii.} For all $X\in\mathcal{M}(U)$ it holds 
\begin{equation}
\tau\underset{g^{-1}}{}(\partial_{o}\wedge X)=\frac{(-1)^{n+1}}{\det [g]}%
\underline{g}(\partial_{o}\lrcorner(\tau X)).  \label{DI.2}
\end{equation}

\textbf{Proof}

It can be deduced by using the identity\footnote{%
In order to prove it we should use the following identities: $I\underset{%
g^{-1}}{}Y=I\underline {g}^{-1}(Y),$ $X\llcorner\underline{t}^{-1}(Y)=%
\underline{t}^{\dagger }(\underline{t}^{*}(X)\llcorner Y)$ and $\underline{t}%
^{-1}(I)=\left. \det\right. ^{-1}[t]I,$ and so forth.} $I\underset{g^{-1}}{}%
Y=\left. \det\right. ^{-1}[g]\underline{g}(IY)$ and Eq.(\ref{DI.1}).$%
\blacksquare$

\subsection{Hodge Duality Identities}

Now, we present two noticeable identities which relate the curl $\partial
_{o}\wedge$ to the contracted divergence $\partial_{o}\lrcorner$ involving
the standard and the metric Hodge extensor fields $\star$ and $\underset{g}{%
\star }$.\vspace{0.1in}

\textbf{i.} For all $X\in\mathcal{M}(U)$ it holds 
\begin{equation}
\star^{-1}(\partial_{o}\wedge(\star X))=-\partial_{o}\lrcorner\widehat{X}.
\label{HDI.1}
\end{equation}

\textbf{Proof}

We will use the duality identity given by Eq.(\ref{DI.1}). By using Eq.(\ref%
{OHD.8a}) and Eq.(\ref{OHD.8}), and recalling the identities $\widetilde{%
\partial_{o}\wedge Y}=\partial_{o}\wedge\overline{Y}$ and $\overline{XY}=%
\overline{Y}$ $\overline{X},$ and the obvious property $\tau\overline{\tau}%
=(-1)^{n},$ we get 
\begin{align*}
\star^{-1}(\partial_{o}\wedge(\star X)) & =\tau\widetilde{(\partial
_{o}\wedge(\widetilde{X}\tau))}=\tau(\partial_{o}\wedge(\overline{\tau }%
\widehat{X})), \\
& =(-1)^{n+1}\partial_{o}\lrcorner(\tau\overline{\tau}\widehat{X}%
)=-\partial_{o}\lrcorner\widehat{X}.\blacksquare
\end{align*}

\textbf{ii.} For all $X\in\mathcal{M}(U)$ it holds 
\begin{equation}
\underset{g}{\star}^{-1}(\partial_{o}\wedge(\underset{g}{\star}X))=-\frac {1%
}{\sqrt{\left| \det[g]\right| }}\underline{g}(\partial_{o}\lrcorner (\sqrt{%
\left| \det[g]\right| }\underline{g}^{-1}(\widehat{X}))).  \label{HDI.2}
\end{equation}

\textbf{Proof}

We will use the duality identity given by Eq.(\ref{DI.2}). A straightforward
calculation using Eq.(\ref{OHD.9a}) and Eq.(\ref{OHD.9}) allows us to get 
\begin{align*}
\underset{g}{\star}^{-1}(\partial_{o}\wedge(\underset{g}{\star}X)) &
=(-1)^{q}\sqrt{\left\vert \det[g]\right\vert }\tau\underset{g^{-1}}{}%
\widetilde{(\partial_{o}\wedge(\widetilde{X}\underset{g^{-1}}{}\underset{g}{%
\tau}))} \\
& =(-1)^{q}\sqrt{\left\vert \det[g]\right\vert }\tau\underset{g^{-1}}{}%
\partial_{o}\wedge(\sqrt{\left\vert \det[g]\right\vert }\overline{\tau }%
\underline{g}^{-1}(\widehat{X})) \\
& =(-1)^{q}\sqrt{\left\vert \det[g]\right\vert }\frac{(-1)^{n+1}}{\det [g]}%
\underline{g}(\partial_{o}\lrcorner(\sqrt{\left\vert \det[g]\right\vert }\tau%
\overline{\tau}\underline{g}^{-1}(\widehat{X}))), \\
& =-\frac{1}{\sqrt{\left\vert \det[g]\right\vert }}\underline{g}(\partial
_{o}\lrcorner(\sqrt{\left\vert \det[g]\right\vert }\underline{g}^{-1}(%
\widehat{X}))).\blacksquare
\end{align*}

\subsection{Ordinary Hodge Coderivative Operators}

We introduce the so-called \emph{standard Hodge derivative operator} $\delta:%
\mathcal{M}(U)\rightarrow\mathcal{M}(U)$ such that 
\begin{equation}
\delta X=\star^{-1}(\partial_{o}\wedge(\star\widehat{X})).  \label{OHO.1}
\end{equation}

Its basic property is 
\begin{equation}
\delta X=-\partial_{o}\lrcorner X.  \label{OHO.1a}
\end{equation}
It immediately follows from the Hodge duality identity given by Eq.(\ref%
{HDI.1}).

We introduce the so-called \emph{metric Hodge coderivative operator} $%
\underset{g}{\delta}:\mathcal{M}(U)\rightarrow\mathcal{M}(U)$ such that 
\begin{equation}
\underset{g}{\delta}X=\text{ }\underset{g}{\star}^{-1}(\partial_{o}\wedge(%
\underset{g}{\star}\widehat{X})).  \label{OHO.2}
\end{equation}

It satisfies the basic property 
\begin{equation}
\underset{g}{\delta}X=-\frac{1}{\sqrt{\left\vert \det[g]\right\vert }}%
\underline{g}(\partial_{o}\lrcorner(\sqrt{\left\vert \det[g]\right\vert }%
\underline{g}^{-1}(X))),  \label{OHO.2a}
\end{equation}
which is an immediate consequence of the Hodge duality identity given by Eq.(%
\ref{HDI.2}).

\section{Levi-Civita Geometric Structure}

We recall \cite{2} that the \emph{Levi-Civita connection field} is the
smooth vector elementary $2$-extensor field on $U,$ namely $\lambda,$
defined by 
\begin{equation}
\lambda(a,b)=\frac{1}{2}g^{-1}\circ(a\cdot\partial_{o}g)(b)+\omega _{0}(a)%
\underset{g}{\times}b,  \label{LGS.1}
\end{equation}
where $\omega_{0}$ is the smooth $(1,2)$-extensor field on $U$ given by 
\begin{equation}
\omega_{0}(a)=-\frac{1}{4}\underline{g}^{-1}(\partial_{b}\wedge\partial
_{c})a\cdot((b\cdot\partial_{o}g)(c)-(c\cdot\partial g_{o})(b)).
\label{LGS.1a}
\end{equation}

Such $\omega_{0}$ satisfies 
\begin{equation}
\omega_{0}(a)\underset{g}{\times}b\underset{g}{\cdot}c=\frac{1}{2}%
a\cdot((b\cdot\partial_{o}g)(c)-(c\cdot\partial_{o}g)(b)).  \label{LGS.1b}
\end{equation}

The \emph{open set} $U$ endowed with $\lambda$ and $g,$ namely $(U,\lambda
,g),$ is a \emph{geometric structure} on $U$, a statement that means that 
\emph{Levi-Civita parallelism structure} $(U,\lambda)$ is compatible with
the metric structure $(U,g).$ Or equivalently, the pair of $a$-\emph{DCDO's }%
associated to $(U,\lambda),$ namely $(D_{a}^{+},D_{a}^{-}),$ is $g$%
-compatible.

The \emph{Levi-Civita} $a$-\emph{DCDO's }$D_{a}^{+}$ and $D_{a}^{-}$ are
defined by 
\begin{align}
D_{a}^{+}X & =a\cdot\partial_{o}X+\Lambda_{a}(X),  \label{LGS.2a} \\
D_{a}^{-}X & =a\cdot\partial_{o}X-\Lambda_{a}^{\dagger}(X).  \label{LGS.2b}
\end{align}
Note that $\Lambda_{a}$ is the so-called \emph{generalized} of $\lambda_{a}.$
The latter is the so-called $a$-\emph{directional connection field}
associated to $\lambda,$ given by $\lambda_{a}(b)=\lambda(a,b).$

We present now two pairs of noticeable properties of $\lambda.\vspace{0.1in} 
$

\textbf{i. }The \emph{scalar divergence} of $\lambda_{a}(b)$ with respect to
the first variable, namely $\partial_{a}\cdot\lambda_{a}(b),$ and the \emph{%
curl} of $\lambda_{a}^{\dagger}(b)$ with respect to the first variable,
namely $\partial_{a}\wedge\lambda_{a}^{\dagger}(b),$ are given by 
\begin{align}
\partial_{a}\cdot\lambda_{a}(b) & =\frac{1}{\sqrt{\left| \det[g]\right| }}%
b\cdot\partial_{o}\sqrt{\left| \det[g]\right| },  \label{LGS.3a} \\
\partial_{a}\wedge\lambda_{a}^{\dagger}(b) & =0.  \label{LGS.3b}
\end{align}

\textbf{Proof}

In order to prove the first result we will use the formula $\tau^{\ast
}(\partial_{n})\cdot(a\cdot\partial_{o}\tau)(n)=\left. \det\right.
^{-1}[\tau]a\cdot\partial_{o}\det[\tau]$, valid for all non-singular smooth $%
(1,1)$-extensor field $\tau.$ By using the symmetry property $%
(g^{-1})^{\dagger}=g^{-1}$ and Eq.(\ref{LGS.1b}), we can write 
\begin{align*}
\partial_{a}\cdot\lambda_{a}(b) & =\frac{1}{2}\partial_{a}\cdot(g^{-1}%
\circ(a\cdot\partial_{o}g)(b)+\partial_{a}\cdot(\omega_{0}(a)\underset{g}{%
\times}b) \\
& =\frac{1}{2}g^{-1}(\partial_{a})\cdot(a\cdot\partial_{o}g)(b)+\omega
_{0}(g^{-1}(\partial_{a}))\underset{g}{\times}b\underset{g}{\cdot}a \\
& =\frac{1}{2}g^{-1}(\partial_{a})\cdot(b\cdot\partial_{o}g)(a)=\frac{1}{2}%
\frac{1}{\det[g]}b\cdot\partial_{o}\det[g].
\end{align*}
Then, recalling the identity $(f)^{-1}b\cdot\partial_{o}f=2(\left\vert
f\right\vert )^{-1/2}b\cdot\partial_{o}(\left\vert f\right\vert )^{1/2}$
valid for all \emph{non-zero} $f\in\mathcal{S}(U)$, the expected result
immediately follows.

To prove the second result we only need to take into account the \emph{%
symmetry property} $\lambda_{a}(b)=\lambda_{b}(a).$ We have 
\begin{equation*}
\partial_{a}\wedge\lambda_{a}^{\dagger}(b)=\partial_{a}\wedge\partial
_{n}(n\cdot\lambda_{a}^{\dagger}(b))=\partial_{a}\wedge\partial_{n}(%
\lambda_{a}(n)\cdot b)=0.\blacksquare
\end{equation*}

\textbf{ii. }The \emph{left contracted divergence} of $\Lambda_{a}(X)$ with
respect to $a,$ namely $\partial_{a}\lrcorner\Lambda_{a}(X),$ and the \emph{%
curl} of $\Lambda_{a}^{\dagger}(X)$ with respect to $a,$ namely $%
\partial_{a}\wedge\Lambda_{a}^{\dagger}(X),$ are given by 
\begin{align}
\partial_{a}\lrcorner\Lambda_{a}(X) & =\frac{1}{\sqrt{\left| \det [g]\right| 
}}(\partial_{o}\sqrt{\left| \det[g]\right| })\lrcorner X,  \label{LGS.4a} \\
\partial_{a}\wedge\Lambda_{a}^{\dagger}(X) & =0.  \label{LGS.4b}
\end{align}

\textbf{Proof}

In order to prove the first ones we will use the multivector identities $%
v\lrcorner(X\wedge Y)=(v\lrcorner X)\wedge Y+\widehat{X}\wedge(v\lrcorner Y)$
and $X\lrcorner(Y\lrcorner Z)=(X\wedge Y)\lrcorner Z,$ where $v\in \mathcal{U%
}_{o}$ and $X,Y,Z\in\bigwedge\mathcal{U}_{o}.$ A straightforward calculation
using Eq.(\ref{LGS.3a}) allows us to get 
\begin{align*}
\partial_{a}\lrcorner\Lambda_{a}(X) &
=\partial_{a}\cdot\lambda_{a}(\partial_{n})(n\lrcorner
X)-\lambda_{a}(\partial_{n})\wedge(\partial _{a}\lrcorner(n\lrcorner X)) \\
& =\frac{1}{\sqrt{\left\vert \det[g]\right\vert }}(\partial_{n}\cdot
\partial_{o}\sqrt{\left\vert \det[g]\right\vert })(n\lrcorner X)-\lambda
_{a}(\partial_{n})\wedge((\partial_{a}\wedge n)\lrcorner X), \\
& =\frac{1}{\sqrt{\left\vert \det[g]\right\vert }}\partial_{n}(n\cdot
\partial_{o}\sqrt{\left\vert \det[g]\right\vert })\lrcorner X-\lambda
_{\partial_{a}}(\partial_{n})\wedge((a\wedge n)\lrcorner X).
\end{align*}
Then, recalling the identity $\partial_{n}(n\cdot\partial_{o}Y)=%
\partial_{o}Y $ and the symmetry property $\lambda_{a}(b)=\lambda_{b}(a),$
we get the required result.

The proof of the second property follows immediately by using Eq.(\ref%
{LGS.3b}), and recalling that the adjoint of generalized is equal to the
generalized of adjoint, i.e., $\Lambda_{a}^{\dagger}(X)=\lambda_{a}^{%
\dagger}(\partial _{n})\wedge(n\lrcorner X)$.$\blacksquare$

\subsection{Levi-Civita Derivatives}

We introduce now the \emph{canonical covariant divergence operator}, $%
D^{+}\lrcorner:\mathcal{M}(U)\rightarrow\mathcal{M}(U)$ such that 
\begin{equation}
D^{+}\lrcorner X=\partial_{a}\lrcorner(D_{a}^{+}X),  \label{LCD.1}
\end{equation}
i.e., $D^{+}\lrcorner
X=e^{\mu}\lrcorner(D_{e_{\mu}}^{+}X)=e_{\mu}\lrcorner(D_{e^{\mu}}^{+}X),$
where $(\{e_{\mu}\},\{e^{\mu}\})$ is any pair of \emph{canonical reciprocal
frame fields} on $U.$

Its basic property is 
\begin{equation}
D^{+}\lrcorner X=\frac{1}{\sqrt{\left\vert \det[g]\right\vert }}\partial
_{o}\lrcorner(\sqrt{\left\vert \det[g]\right\vert }X).  \label{LCD.1a}
\end{equation}
Indeed, using the identity $\partial_{a}\lrcorner(a\cdot\partial
_{o}X)=\partial_{o}\lrcorner X$ and Eq.(\ref{LGS.4a}) into the definition
given by Eq.(\ref{LGS.2a}), we get 
\begin{align*}
D^{+}\lrcorner X & =\partial_{a}\lrcorner(a\cdot\partial_{o}X)+\partial
_{a}\lrcorner\Lambda_{a}(X), \\
& =\partial\lrcorner X+\frac{1}{\sqrt{\left\vert \det[g]\right\vert }}%
(\partial_{o}\sqrt{\left\vert \det[g]\right\vert })\lrcorner X.
\end{align*}
So, by recalling the identity $\partial_{o}\lrcorner(fY)=(\partial
_{o}f)\lrcorner X+f(\partial_{o}\lrcorner X),$ for all $f\in\mathcal{S}(U)$
and $Y\in\mathcal{M}(U),$ we can get the proof for this remarkable property.

The so-called \emph{metric covariant divergence}, \emph{covariant curl }and 
\emph{metric covariant gradient operators}, namely $D^{-}\underset{g^{-1}}{%
\lrcorner},$ $D^{-}\wedge$ and $D^{-}\underset{g^{-1}}{},$ map also smooth
multivector fields to smooth multivector fields and are defined by 
\begin{align}
D^{-}\underset{g^{-1}}{\lrcorner}X & =\partial_{a}\underset{g^{-1}}{\lrcorner%
}(D_{a}^{-}X)=g^{-1}(\partial_{a})\lrcorner(D_{a}^{-}X),  \label{LCD.2a} \\
D^{-}\wedge X & =\partial_{a}\wedge(D_{a}^{-}X),  \label{LCD.2b} \\
D^{-}\underset{g^{-1}}{}X & =\partial_{a}\underset{g^{-1}}{}(D_{a}^{-}X).
\label{LCD.2c}
\end{align}
The relationship among these operators is given by 
\begin{equation}
D^{-}\underset{g^{-1}}{}X=D^{-}\underset{g^{-1}}{\lrcorner}X+D^{-}\wedge X.
\label{LCD.3}
\end{equation}

The basic properties of the metric covariant divergence are 
\begin{align}
D^{-}\underset{g^{-1}}{\lrcorner}X & =\underline{g}(D^{+}\lrcorner 
\underline{g}^{-1}(X)),  \label{LCD.4a} \\
D^{-}\underset{g^{-1}}{\lrcorner}X & =\frac{1}{\sqrt{\left| \det[g]\right| }}%
\underline{g}(\partial_{o}\lrcorner(\sqrt{\left| \det[g]\right| }\underline{g%
}^{-1}(X))).  \label{LCD.4b}
\end{align}

Eq.(\ref{LCD.4a}) follows from the fundamental property $D_{a}^{-}\underline{%
g}(X)=\underline{g}(D_{a}^{+}X)$ which holds for any $g$-compatible pair of $%
a$-\emph{DCDO's, }by using the identity $X\lrcorner\underline {t}(Y)=%
\underline{t}(\underline{t}^{\dagger}(X)\lrcorner Y).$\emph{\ }Eq.(\ref%
{LCD.4b}) is deduced by using Eq.(\ref{LCD.4a}) and Eq.(\ref{LCD.1a}).

A remarkable property which follows from Eq.(\ref{LCD.4b}) is 
\begin{equation}
D^{-}\underset{g^{-1}}{\lrcorner}(D^{-}\underset{g^{-1}}{\lrcorner}X)=0.
\label{LCD.4b1}
\end{equation}
In order to prove it we should use the known identity $\partial_{o}%
\lrcorner(\partial_{o}\lrcorner Y)=0.$

By comparing Eq.(\ref{OHO.2a}) and Eq.(\ref{LCD.4b}) we get 
\begin{equation}
D^{-}\underset{g^{-1}}{\lrcorner}X=-\underset{g}{\delta}X.  \label{LCD.4b2}
\end{equation}

The basic property for the covariant curl is 
\begin{equation}
D^{-}\wedge X=\partial_{o}\wedge X.  \label{LCD.5}
\end{equation}
It immediately follows by using the identity $\partial_{a}\wedge
(a\cdot\partial_{o}X)=\partial_{o}\wedge X$ and Eq.(\ref{LGS.4b}) into the
definition given by Eq.(\ref{LGS.2b}).

The four \emph{derivative-like operators} defined by Eq.(\ref{LCD.1}), Eq.(%
\ref{LCD.2a}), Eq.(\ref{LCD.2b}) and Eq.(\ref{LCD.2c}) will be called the 
\emph{Levi-Civita derivatives. }They are involved in three useful identities
which are used in the \emph{Lagrangian theory of multivector fields}. These
are, 
\begin{align}
& (\partial_{o}\wedge X)\underset{g^{-1}}{\cdot}Y+X\underset{g^{-1}}{\cdot }%
(D^{-}\underset{g^{-1}}{\lrcorner}Y)  \notag \\
& =\frac{1}{\sqrt{\left\vert \det[g]\right\vert }}\partial_{o}\cdot (\sqrt{%
\left\vert \det[g]\right\vert }\partial_{n}(n\wedge X)\underset{g^{-1}}{\cdot%
}Y),  \label{LCD.6a} \\
& (D^{-}\underset{g^{-1}}{\lrcorner}X)\underset{g^{-1}}{\cdot}Y+X\underset{%
g^{-1}}{\cdot}(\partial_{o}\wedge Y)  \notag \\
& =\frac{1}{\sqrt{\left\vert \det[g]\right\vert }}\partial_{o}\cdot (\sqrt{%
\left\vert \det[g]\right\vert }\partial_{n}(n\underset{g^{-1}}{\lrcorner}X)%
\underset{g^{-1}}{\cdot}Y),  \label{LCD.6b} \\
& (D^{-}\underset{g^{-1}}{}X)\underset{g^{-1}}{\cdot}Y+X\underset{g^{-1}}{%
\cdot}(D^{-}\underset{g^{-1}}{}Y)  \notag \\
& =\frac{1}{\sqrt{\left\vert \det[g]\right\vert }}\partial_{o}\cdot (\sqrt{%
\left\vert \det[g]\right\vert }\partial_{n}(n\underset{g^{-1}}{}X)\underset{%
g^{-1}}{\cdot}Y).  \label{LCD.6c}
\end{align}

\subsection{Gauge Derivatives}

Let $h$ be a gauge metric field for $g$ \cite{3}$.$ This statement means
that there is a smooth $(1,1)$-extensor field $h$ such that $%
g=h^{\dagger}\circ \eta\circ h,$ where $\eta$ is an orthogonal metric field
with the same signature as $g.$ As we know, associated to a Levi-Civita $a$-%
\emph{DCDO's,} namely $(D_{a}^{+},D_{a}^{-}),$ there must be an unique $\eta$%
-compatible pair of $a$-\emph{DCDO's, }namely\emph{\ }$($\DJ $_{a}^{+},$\DJ $%
_{a}^{-}),$ given by the following formulas 
\begin{align}
\text{\DJ }_{ha}^{+}X & =\underline{h}(D_{a}^{+}\underline{h}^{-1}(X)),
\label{GD.1a} \\
\text{\DJ }_{h^{\ast}a}^{-}X & =\underline{h}^{\ast}(D_{a}^{+}\underline {h}%
^{\dagger}(X)).  \label{GD.1b}
\end{align}
These equations say that $($\DJ $_{ha}^{+},$\DJ $_{h^{\ast}a}^{-})$ is the $%
h $-\emph{deformation} of $(D_{a}^{+},D_{a}^{-}).$ So, \DJ $_{ha}^{+}$ and 
\DJ\ $_{h^{\ast}a}^{-}$ will be called the \emph{gauge covariant derivatives}
associated to $D_{a}^{+}$ and $D_{a}^{-}.$

We present here two noticeable properties of \DJ $_{ha}^{+}.\vspace{0.1in}$

\textbf{i.} For all $a,b,c\in\mathcal{V}(U),$ it holds 
\begin{equation}
(\text{\DJ }_{ha}^{+}b)\underset{\eta}{\cdot}c=[h(a),h^{-1}(b),h^{-1}(c)].
\label{GD.2}
\end{equation}

\textbf{Proof}

By using Eq.(\ref{GD.1a}), the identity $(D_{a}^{+}b)\cdot c=\QDATOPD{\{}{\}%
}{c}{a,b},$ and the definition of the Christoffel operator of second kind,
i.e., $\QDATOPD{\{}{\}}{c}{a,b}=[a,b,g^{-1}(c)],$ we can write 
\begin{equation*}
(\text{\DJ }_{ha}^{+}b)\underset{\eta}{\cdot}c=D_{a}^{+}(h^{-1}(b))\cdot
h^{\dagger}\circ\eta(c)=\QATOPD{\{}{\}}{h^{\dagger}\circ\eta(c)}{%
h(a),h^{-1}(b)}=[h(a),h^{-1}(b),g^{-1}\circ h^{\dagger}\circ\eta(c)],
\end{equation*}
hence, by recalling that $g^{-1}=h^{-1}\circ\eta\circ h^{\ast},$ the
required result immediately follows.$\blacksquare$

\textbf{ii.} There exists a smooth $(1,2)$-extensor field on $U_{o},$ namely 
$\Omega_{0},$ such that 
\begin{equation}
\text{\DJ }_{a}^{+}X=a\cdot\partial_{o}X+\Omega_{0}(a)\underset{\eta}{\times 
}X.  \label{GD.3}
\end{equation}
Such $\Omega_{0}$ is given by 
\begin{equation}
\Omega_{0}(a)=-\frac{1}{2}\underline{\eta}(\partial_{b}\wedge\partial
_{c})[a,h^{-1}(b),h^{-1}(c)].  \label{GD.4}
\end{equation}

\textbf{Proof}

Firstly, we must prove a particular case of the above property, i.e., 
\begin{equation}
\text{\DJ }_{a}^{+}b=a\cdot\partial_{o}b+\Omega_{0}(a)\underset{\eta}{\times 
}b.  \label{GD.3a}
\end{equation}

We will use the following properties of the Christoffel operator of first
kind $[a,b,c]-[a,c,b]=2[a,b,c]-a\cdot\partial_{o}(b\underset{g}{\cdot}c),$
and $[a,b+b^{\prime},c]=[a,b,c]+[a,b^{\prime},c]$ and $[a,fb,c]=f[a,b,c]+(a%
\cdot \partial_{o}f)b\underset{g}{\cdot}c.$ A straightforward calculation
allows us to get 
\begin{align*}
\Omega_{0}(a)\underset{\eta}{\times}b\underset{\eta}{\cdot}c & =(b\wedge c)%
\underset{\eta}{\cdot}\Omega_{0}(a)=\frac{1}{2}(b\wedge c)\cdot(\partial
_{p}\wedge\partial_{q})[a,h^{-1}(p),h^{-1}(q)] \\
& =\frac{1}{2}\det\left[ 
\begin{array}{cc}
b\cdot\partial_{p} & b\cdot\partial_{q} \\ 
c\cdot\partial_{p} & c\cdot\partial_{q}%
\end{array}
\right] [\ldots] \\
& =\frac{1}{2}b\cdot\partial_{p}c\cdot\partial_{q}[\ldots]-b\cdot\partial
_{q}c\cdot\partial_{p}[\ldots] \\
& =\frac{1}{2}b\cdot\partial_{p}c\cdot%
\partial_{q}([a,h^{-1}(p),h^{-1}(q)]-[a,h^{-1}(q),h^{-1}(p)]) \\
& =\frac{1}{2}b\cdot\partial_{p}c\cdot\partial_{q}(2[\ldots]-a\cdot
\partial_{o}(h^{-1}(p)\underset{g}{\cdot}h^{-1}(q))) \\
& =b\cdot\partial_{p}c\cdot\partial_{q}(p\cdot b_{\mu}[a,h^{-1}(b_{\mu
}),h^{-1}(c)]+(a\cdot\partial_{o}p)\underset{\eta}{\cdot}c) \\
& =b\cdot b_{\mu}[a,h^{-1}(b_{\mu}),h^{-1}(c)]+(a\cdot\partial_{o}b)\underset%
{\eta}{\cdot}c-(a\cdot\partial_{o}b)\underset{\eta}{\cdot}c \\
& =[a,h^{-1}(b),h^{-1}(c)]-(a\cdot\partial_{o}b)\underset{\eta}{\cdot}c,
\end{align*}
hence, by using Eq.(\ref{GD.2}), the particular case given by Eq.(\ref{GD.3a}%
) immediately follows. Now, we can indeed prove the general case of the
above property .

As we can see, the $a$-directional connection field for $($\DJ $_{a}^{+},$%
\DJ $_{a}^{-})$ is given by $b\mapsto\Omega_{0}(a)\underset{\eta}{\times}b.$
Then, its generalized (extensor field) must be given by $X\mapsto(\Omega
_{0}(a)\underset{\eta}{\times}\partial_{b})\wedge(b\lrcorner X).$ But, by
recalling the noticeable identity $(B\underset{\eta}{\times}\partial
_{b})\wedge(b\lrcorner X)=B\underset{\eta}{\times}X,$ where $B\in\bigwedge
^{2}\mathcal{U}_{o}$ and $X\in\bigwedge\mathcal{U}_{o},$ we find that it can
be written as $X\mapsto\Omega_{0}(a)\underset{\eta}{\times}X.\blacksquare$

We introduce now the \emph{gauge covariant divergence}, \emph{gauge
covariant curl} and \emph{gauge covariant gradient operators}, namely \DJ $%
^{-}\underset{\eta}{\lrcorner},$ \DJ $^{-}\wedge$ and \DJ $^{-}\underset{\eta%
}{}.$ They all map smooth multivector fields to smooth multivector fields
and are defined by 
\begin{align}
\text{\DJ }^{-}\underset{\eta}{\lrcorner}X & =h^{\ast}(\partial _{a})%
\underset{\eta}{\lrcorner}(\text{\DJ }_{h^{\ast}a}^{-}X),  \label{GD.5a} \\
\text{\DJ }^{-}\wedge X & =h^{\ast}(\partial_{a})\wedge(\text{\DJ }%
_{h^{\ast}a}^{-}X),  \label{GD.5b} \\
\text{\DJ }^{-}\underset{\eta}{}X & =h^{\ast}(\partial_{a})\underset{\eta}{}(%
\text{\DJ }_{h^{\ast}a}^{-}X).  \label{GD.5c}
\end{align}
It is obvious that the relationship among them is given by 
\begin{equation}
\text{\DJ }^{-}\underset{\eta}{}X=\text{\DJ }^{-}\underset{\eta}{\lrcorner }%
X+\text{\DJ }^{-}\wedge X.  \label{GD.6}
\end{equation}

Their basic properties are given by the following \emph{golden formulas} 
\begin{align}
\underline{h}^{\ast}(D^{-}\underset{g^{-1}}{\lrcorner}X) & =\text{\DJ }^{-}%
\underset{\eta}{\lrcorner}\underline{h}^{\ast}(X),  \label{GD.7a} \\
\underline{h}^{\ast}(D^{-}\wedge X) & =\text{\DJ }^{-}\wedge\underline {h}%
^{\ast}(X),  \label{GD.7b} \\
\underline{h}^{\ast}(D^{-}\underset{g^{-1}}{}X) & =\text{\DJ }^{-}\underset{%
\eta}{}\underline{h}^{\ast}(X).  \label{GD.7c}
\end{align}
These formulas can be proved by using the golden formula deduced in \cite{3}%
, $\underline{h}^{\ast}(X\underset{g^{-1}}{\ast}Y)=\underline{h}^{\ast }(X)%
\underset{\eta}{\ast}\underline{h}^{\ast}(Y),$ where $X,Y\in \bigwedge%
\mathcal{U}_{o},$ and $\underset{g^{-1}}{\ast}$ means either exterior
product or any $g^{-1}$-product of smooth multivector fields, and
analogously for $\underset{\eta}{\ast}$. It is also necessary to take into
account the master formulas $g=h^{\dagger}\circ\eta\circ h$ and $%
g^{-1}=h^{-1}\circ \eta\circ h^{\ast},$ and the relationship between \DJ $%
_{ha}^{+}$ and \DJ $_{h^{\ast}a}^{-},$ i.e., $\underline{\eta}($\DJ $%
_{ha}^{+}\underline {\eta}(X))=$ \DJ $_{h^{\ast}a}^{-}X.$

>From Eq.(\ref{LCD.4b}) and Eq.(\ref{LCD.5}) by using Eq.(\ref{GD.7a}) and
Eq.(\ref{GD.7b}) we find the interesting identities 
\begin{align}
\text{\DJ }^{-}\underset{\eta}{\lrcorner}X & =\frac{1}{\det[h]}\underline{%
\eta\circ h}(\partial_{o}\lrcorner(\det[h]\underline{h^{-1}\circ\eta}(X))),
\label{GD.8a} \\
\text{\DJ }^{-}\wedge X & =\underline{h}^{\ast}(\partial_{o}\wedge 
\underline{h}^{\dagger}(X)).  \label{GD.8b}
\end{align}

\section{Covariant Hodge Coderivative}

Let $(U,\gamma,g)$ be a \emph{geometric structure} on $U,$ and let us denote
by $(\mathcal{D}_{a}^{+},\mathcal{D}_{a}^{-})$ the $g$-\emph{compatible}
pair of $a$-\emph{DCDO's }associated to $(U,\gamma,g).$ We will present two
noticeable properties which are satisfied by the second $a$-\emph{DCDO.}

\textbf{i.} $\underset{g}{\tau}$ is a \emph{covariantly constant} smooth
pseudoscalar field on $U,$ i.e., 
\begin{equation}
\mathcal{D}_{a}^{-}\underset{g}{\tau}=0.  \label{CHC.1}
\end{equation}

\textbf{Proof}

We will use the basic property and the expansion formula given by Eq.(\ref%
{OHD.7a}) and Eq.(\ref{OHD.7b}). By using the Ricci-like theorem for $%
\mathcal{D}_{a}^{-},$ we have that 
\begin{equation*}
(\mathcal{D}_{a}^{-}\underset{g}{\tau})\underset{g^{-1}}{\cdot}\underset{g}{%
\tau}+\underset{g}{\tau}\underset{g^{-1}}{\cdot}(\mathcal{D}_{a}^{-}\underset%
{g}{\tau})=0,
\end{equation*}
i.e., $(\mathcal{D}_{a}^{-}\underset{g}{\tau})\underset{g^{-1}}{\cdot }%
\underset{g}{\tau}=0.$ Then, $\mathcal{D}_{a}^{-}\underset{g}{\tau}%
=(-1)^{q}((\mathcal{D}_{a}^{-}\underset{g}{\tau})\underset{g^{-1}}{\cdot }%
\underset{g}{\tau})\underset{g}{\tau}=0$.$\blacksquare$

\textbf{ii.} For all $X\in\mathcal{M}(U)$ it holds 
\begin{equation}
\mathcal{D}_{a}^{-}(\underset{g}{\tau}\underset{g^{-1}}{*}X)=\underset{g}{%
\tau}\underset{g^{-1}}{*}\mathcal{D}_{a}^{-}(X),  \label{CHC.2}
\end{equation}
where $\underset{g^{-1}}{*}$ means either exterior product or any $g^{-1}$%
-product of smooth multivector fields.

\textbf{Proof}

It can be deduced by using Eq.(\ref{CHC.1}) and the general Leibnitz rule
for $\mathcal{D}_{a}^{-}.\blacksquare$

We introduce the so-called \emph{metric covariant divergence}, \emph{%
covariant curl} and \emph{metric covariant gradient operators}, namely $%
\mathcal{D}^{-}\underset{g^{-1}}{\lrcorner},$ $\mathcal{D}^{-}\wedge$ and $%
\mathcal{D}^{-}\underset{g^{-1}}{},$ which map smooth multivector fields to
smooth multivector fields. They are defined by 
\begin{align}
\mathcal{D}^{-}\underset{g^{-1}}{\lrcorner}X & =\partial_{a}\underset{g^{-1}}%
{\lrcorner}(\mathcal{D}_{a}^{-}X),  \label{CHC.3a} \\
\mathcal{D}^{-}\wedge X & =\partial_{a}\wedge(\mathcal{D}_{a}^{-}X),
\label{CHC.3b} \\
\mathcal{D}^{-}\underset{g^{-1}}{}X & =\partial_{a}\underset{g^{-1}}{}(%
\mathcal{D}_{a}^{-}X).  \label{CHC.3c}
\end{align}
The relationship among them is given by 
\begin{equation}
\mathcal{D}^{-}\underset{g^{-1}}{}X=\mathcal{D}^{-}\underset{g^{-1}}{%
\lrcorner}X+\mathcal{D}^{-}\wedge X.  \label{CHC.3d}
\end{equation}

Now, we will present two noticeable duality identities between $\mathcal{D}%
^{-}\wedge$ and $\mathcal{D}^{-}\underset{g^{-1}}{\lrcorner}.$ One of them
involves the metric Hodge extensor field $\underset{g}{\star}.$

\textbf{iii. }For all $X\in\mathcal{M}(U)$ it holds 
\begin{equation}
\underset{g}{\tau}\underset{g^{-1}}{}(\mathcal{D}^{-}\wedge X)=(-1)^{n+1}%
\mathcal{D}^{-}\underset{g^{-1}}{\lrcorner}(\underset{g}{\tau}\underset{%
g^{-1}}{}X).  \label{CHC.4}
\end{equation}

\textbf{Proof}

We will use the $g^{-1}$-duality identity $I\underset{g^{-1}}{}(a\wedge
Y)=(-1)^{n+1}a\underset{g^{-1}}{\lrcorner}(I\underset{g^{-1}}{}Y),$ where $%
a\in\mathcal{U}_{o},$ $I\in\bigwedge^{n}\mathcal{U}_{o}$ and $Y\in \bigwedge%
\mathcal{U}_{o}.$ A straightforward calculation by taking into account Eq.(%
\ref{CHC.2}) gives 
\begin{align*}
\underset{g}{\tau}\underset{g^{-1}}{}(\mathcal{D}^{-}\wedge X) &
=(-1)^{n+1}\partial_{a}\underset{g^{-1}}{\lrcorner}(\underset{g}{\tau }%
\underset{g^{-1}}{}(\mathcal{D}_{a}^{-}X)) \\
& =(-1)^{n+1}\partial_{a}\underset{g^{-1}}{\lrcorner}(\mathcal{D}_{a}^{-}(%
\underset{g}{\tau}\underset{g^{-1}}{}X)), \\
& =(-1)^{n+1}\mathcal{D}^{-}\underset{g^{-1}}{\lrcorner}(\underset{g}{\tau }%
\underset{g^{-1}}{}X).\blacksquare
\end{align*}

\textbf{iv. }For all $X\in\mathcal{M}(U)$ it holds 
\begin{equation}
\underset{g}{\star}^{-1}(\mathcal{D}^{-}\wedge(\underset{g}{\star }X))=-%
\mathcal{D}^{-}\underset{g^{-1}}{\lrcorner}\widehat{X}.  \label{CHC.5}
\end{equation}

\textbf{Proof}

We will use the duality identity given by Eq.(\ref{CHC.4}). By using Eq.(\ref%
{OHD.9a}) Eq.(\ref{OHD.9}), and recalling the identities $\widetilde{%
\mathcal{D}^{-}\wedge X}=\mathcal{D}^{-}\wedge\overline{X}$ and $\overline{XY%
}=\overline{Y}$ $\overline{X},$ and the obvious property $\underset{g}{\tau}%
\underset{g^{-1}}{}\overline{\underset{g}{\tau}}=(-1)^{n+q},$ we get 
\begin{align*}
\underset{g}{\star}^{-1}(\mathcal{D}^{-}\wedge(\underset{g}{\star}X)) &
=(-1)^{q}\underset{g}{\tau}\underset{g^{-1}}{}\widetilde{(\mathcal{D}%
^{-}\wedge(\widetilde{X}\underset{g^{-1}}{}\underset{g}{\tau}))} \\
& =(-1)^{q}\underset{g}{\tau}\underset{g^{-1}}{}(\mathcal{D}^{-}\wedge(%
\overline{\underset{g}{\tau}}\underset{g^{-1}}{}\widehat{X}) \\
& =(-1)^{q}(-1)^{n+1}\mathcal{D}^{-}\underset{g^{-1}}{\lrcorner}(\underset{g}%
{\tau}\underset{g^{-1}}{}\overline{\underset{g}{\tau}}\underset{g^{-1}}{}%
\widehat{X}), \\
& =-\mathcal{D}^{-}\underset{g^{-1}}{\lrcorner}\widehat{X}.\blacksquare
\end{align*}

We introduce now the \emph{covariant Hodge coderivative operator}, namely $%
\underset{g}{\Delta},$ as defined by $\underset{g}{\Delta}:\mathcal{M}%
(U)\rightarrow\mathcal{M}(U)$ such that 
\begin{equation}
\underset{g}{\Delta}X=\underset{g}{\text{ }\star}^{-1}(\mathcal{D}^{-}\wedge(%
\underset{g}{\star}\widehat{X})).  \label{CHC.6}
\end{equation}
Its basic property is 
\begin{equation}
\mathcal{D}^{-}\underset{g^{-1}}{\lrcorner}X=-\underset{g}{\Delta}X,
\label{CHC.6a}
\end{equation}
which follows trivially from Eq.(\ref{CHC.5}).

\section{Conclusions}

In this paper (the fourth paper in a series of five) we continued the
development of our theory of multivector and extensor calculus on smooth
manifolds. We introduced first the concepts of ordinary Hodge coderivatives,
duality identities, and Hodge coderivative identities. Then, we recalled the
concept of a Levi-Civita geometric structure and the concepts of Levi-Civita
and gauge derivatives. We also introduced the new concept of \emph{covariant}
Hodge coderivative. We studied in details how all these important concepts
are related and obtained several formulas which play a key role in the
Lagrangian theory of multivector and extensor fields on arbitrary
manifolds.We believe that our results clarify many misconceptions in earlier
literature on the subject.\medskip

\textbf{Acknowledgments:} V. V. Fern\'{a}ndez and A. M. Moya are very
grateful to Mrs. Rosa I. Fern\'{a}ndez who gave to them material and
spiritual support at the starting time of their research work. This paper
could not have been written without her inestimable help.

\end{document}